\documentclass[11pt]{article}
\usepackage{latexsym}
\usepackage{amsmath}
\usepackage{amssymb}
\title{\LARGE Functionals depending on curvatures  
with constraints}
\vspace{3cm}
\author{\hphantom{AAAA}\\
\\
\\
{\bf \large Maria Giovanna Mora} \vspace{.3cm}\\
{\bf \large Massimiliano Morini} \vspace{.7cm}\\
{\normalsize S.I.S.S.A.}\vspace{.2cm}\\
{\normalsize via Beirut 2/4, 34014 TRIESTE (Italy)}\vspace{.2cm}\\
{\normalsize e-mail: \texttt{mora@sissa.it,
morini@sissa.it}}
}
\date{}

\makeatletter
\@addtoreset{equation}{section}
\makeatother

\newcommand{\enne}{{\mathbb N}}

\newcommand{\reali}{{\mathbb R}}
\newcommand{\rdue}{{\mathbb R}^{2}}
\newcommand{\rn}{{\mathbb R}^n}
\newcommand{\rtre}{{\mathbb R}^{3}}
\newcommand{\matr}{{\mathbb M}\,^{n\times n}}
\newcommand{\sobdue}{W^{2,\infty}}
\newcommand{\dist}{\hbox{dist}}
\newcommand{\qed}{\hfill $\square$ \bigskip}

\newcommand{\graph}{{\rm graph}\,}
\newcommand{\kconv}{\overset{{\cal K}}{\to}}
\newcommand{\inter}[1]{{\rm Int}\, #1 }

\newcommand{\diam}{{\rm diam}}
\newcommand{\pbar}{\overline{p}}
\newcommand{\xbar}{\overline{x}}
\newcommand{\ybar}{\overline{y}}
\newcommand{\esse}[2]{{\cal S}_{#2}^{#1}}
\newcommand{\proj}[2]{\Pi_{#1}^{#2}}
\newcommand{\avec}[2]{a_{#1}^{#2}}
\newcommand{\bvec}[2]{b_{#1}^{#2}}

\newtheorem{proposition}{Proposition}[section]
\newtheorem{lemma}[proposition]{Lemma}
\newtheorem{theorem}[proposition]{Theorem}

\newtheorem{definition}[proposition]{Definition}
\newtheorem{example}[proposition]{Example}

\begin{document}
\maketitle

\vspace{2.5cm}
\begin{center}
{\bf Abstract}
\end{center}
\vspace{.1cm}
\begin{center}
\begin{minipage}{12.5cm}
\noindent {\normalsize
We deal with a family of functionals depending on curvatures
and we
prove for them compactness and semicontinuity properties
in the class of closed and bounded sets 
which satisfy a uniform exterior and interior sphere condition.
We apply the results
to state an existence theorem for the Nitzberg and
Mumford problem under this additional constraint.
}

\thispagestyle{empty}

\end{minipage}
\end{center}

\begin{center}
\vfill
{\large Ref. S.I.S.S.A. 85/99/M  (July 99)}
\end{center}
\pagebreak
\setcounter{page}{1}

\section{Introduction}

In this paper we are dealing with geometrical functionals of the form
\begin{equation}\label{effe}
F(E)= \int_{\partial E} \varphi(K_1, \ldots, K_{n-1}) \,d {\cal H}^{n-1},
\end{equation}
where $\varphi: \reali^{n-1} \to \reali$ is a given convex function,
$E$ varies in a class of sufficiently regular closed subsets of $\rn$, 
$K_1, \ldots, K_{n-1}$ denote the elementary symmetric curvatures 
of $\partial E$ (see (\ref{symm})), and ${\cal H}^{n-1}$ is the
$(n-1)$-dimensional Hausdorff measure.

In \cite{bel-dm-pao} {\sc Bellettini}, {\sc Dal Maso}, {\sc Paolini} 
studied the functional $F$
in the case $n=2$ and $\varphi(\kappa) = 1 + |\kappa |^p$, 
where $\kappa$ denotes the curvature of $\partial E$,
and remarked that $F$
does not have the right compactness properties in its natural class of
definition, composed of all closed sets $E$ whose boundary is of class 
$W^{2,p}$: simple examples show that there exist sets of class $C^{\infty}$,
except for a finite number of cusp points (the functional then 
is not naturally defined on them), 
which can be approximated by a sequence of sets of class $C^{\infty}$, 
whose boundaries have bounded curvature. Moreover, they show
that the lower semicontinuous envelope, $\overline{F}$, of $F$
with respect to the $L^1$-topology cannot be represented as an integral 
of the form
$$\int_{\partial E} f(\kappa) \, d{\cal H}^1,$$
and that the fact that a set $E$ belongs to the domain of $\overline{F}$
depends on the global structure of $E$. 
For instance, if $\partial E$ is smooth
except for a finite number $k$ of cusp points, 
then $\overline{F}(E)<+\infty$ if and
only if $k$ is even.

The idea of this work is to modify the domain of $F$ 
by introducing some suitable 
constraints.

Fixed $R>0$, we choose as domain of $F$ the class
\begin{multline}\label{defi}
{\cal U}_{R}  =   \left\{  \right.   E \subset \rn, E \ \hbox{closed
and bounded} : 
 \forall p \in \partial E \  \exists \, p', p'':  \\
p \in \partial B(p', R) \cap \partial B(p'', R),
\left.  B(p', R) \subset E ,\, B(p'', R)\cap E= \emptyset \right\},
\end{multline}
where $B(q, R)$ denotes the open ball
centred at $q$ of radius $R$; we will say, equivalently, 
that ${\cal U}_{R}$ is the class of all closed and bounded subsets of $\rn$, which satisfy 
{\it the exterior and interior sphere condition} 
with radius $R$ at every point of the boundary.
Note that the introduced constraint has a nonlocal effect on the thickness, 
which cannot be too small, and a local effect on the curvatures,
which are bounded from above by a constant depending only on $R$.
Remark also that this upper bound on the curvatures goes to infinity, when
$R$ tends to $0$.

In the class ${\cal U}_{R}$ the pathological phenomena 
described above cannot occur; indeed, they are
related to the existence of approximating sequences of sets having regions 
with vanishing thickness or
different connected components whose distance goes to $0$.

In Section 2 we study the regularity of sets belonging to ${\cal U}_R$,
showing that the functional in (\ref{effe}) is well defined.
In Sections 3 and 4 we prove compactness and semicontinuity results for $F$
in ${\cal U}_R$.
In Section 5 we consider the case $n=2$ and
we apply the theorems of Sections 3 and 4 to 
show the existence of a solution to the variational problem
\begin{multline}\label{min2dim}
\min\displaystyle \left\{ 
\sum_{i=1}^{k}\left( 
\alpha \int_{E'_{i} \cap \Omega} |g - 
g_{E'_{i} \cap \Omega}|^{2} dx 
  + \beta {\cal L}^{2}(E_{i}) + \gamma \int_{\partial E_{i}} 
\varphi ( \kappa ) 
d{\cal H}^{1}\right) \right.+ \\ 
\left. \hphantom{\sum \int}
+\alpha \int_{\Omega \setminus \cup_{i=1}^{k} E_{i}} 
|g - g_{\Omega 
\setminus \cup_{i=1}^{k} E_{i}}|^{2} dx : 
 \, E_1, \ldots, E_k \in {\cal U}_{R} \right\},
\end{multline}
where $\Omega$ is a bounded subset of $\rdue$, $\alpha, \beta,
\gamma$ are positive parameters,
$E'_i:=E_i\setminus \cup_{j=1}^{i-1}E_j$, $g$ is a function in 
$L^{2}(\Omega )$.
This functional was proposed by {\sc Nitzberg} and {\sc Mumford} 
as a variant of the
{\sc Mumford} and {\sc Shah} image segmentation model, 
allowing regions to overlap (for
further information about this model, see \cite{nit-mum-shi}).
In this framework the constant $R$ can be interpreted 
as a resolution parameter
of the segmented image: the thickness of the reconstructed 
objects has to be greater or equal to $2R$.
We conclude the section by giving an example 
of non trivial minimizer for a functional of the 
form as in (\ref{min2dim}).

\section{Preliminary results}

In this section we investigate the regularity of sets belonging
to the class ${\cal U}_R$ introduced in (\ref{defi}) and we show that
the functional (\ref{effe}) is well defined in this class.

Let us fix first some notation. If $E$ belongs to ${\cal U}_R$
and $p\in \partial E$, we denote the centres of the interior 
and exterior
balls associated to $p$ by $p'$ and $p''$ respectively, as
in (\ref{defi}); moreover, we call
$\esse{p}{E}$ the
class of all coordinate systems centred at $p$ 
such that the vector $\frac{1}{2R} (p'' - p')$
coincides with  
the $n$-th vector of the coordinate basis.

\begin{proposition}\label{regular}
There exists a constant $\rho > 0$ (depending only on $R$), such that
for every $E\in {\cal U}_{R}$ and 
for every $p_0 \in \partial E$, if we call
$C$ the cylinder $\{ x\in \reali^{n-1} :|x| < \rho  \} { \times } 
] { - } R,R[$ expressed with respect to a coordinate
system belonging to $\esse{p_0}{E}$,
then $\partial E \cap C$ is 
the subgraph of a function $f$ belonging to 
$\sobdue ( \{ x\in \reali^{n-1} :|x| < \rho  \}  )$. 
Moreover, the $\sobdue$-norm of $f$ is bounded by a constant depending 
only on $R$ (independent on $p_0$, on $E$ and on the choice of  the coordinate system in $\esse{p_0}{E}$).
\end{proposition}

{\sc Proof.} We first perform the proof in the case $n=2$
showing that $\rho = \sqrt{3}R/2$ is a good choice.

Let $E$ be in ${\cal U}_R$ and let $p_0$ belong to $\partial E$.
Let us consider 
a coordinate system belonging to $\esse{p_0}{E}$.
We can reduce to work 
in the cylinder $C^+ :=[0, \sqrt{3}R/2 [ \times ]{ - }R, R[$.
For the proof we need the following lemma.
 
\begin{lemma} \label{grafico}
Let $\pbar =(\xbar, \ybar)$ be in 
$\partial E \cap C^+$. If we call $\alpha(\pbar)$
the angle in $[0, \pi[$ between
the $x$-axis and the tangent line 
to $B(\pbar', R)$ and $B(\pbar'',R)$ at $\pbar$,
then
\begin{equation}\label{cos}
|\cos \alpha(\pbar) | \geq  \frac{\sqrt{R^2 - \xbar^2}}{R}.
\end{equation}
Moreover, either $B(\pbar', R)$ or $B(\pbar'', R)$ 
contains the whole segment $\{ \xbar \} \times 
] \ybar, R - \sqrt{R^2 - \xbar^2}].$
\end{lemma}

{\sc Proof.} Let us suppose 
by contradiction that (\ref{cos}) does not hold; hence,
\begin{equation} \label{sin}
\sin \alpha(\pbar) > \frac{\xbar}{R}.
\end{equation}
The point $q:=(\xbar - R \sin \alpha(\pbar), \ybar + R \cos \alpha
(\pbar))$ 
must coincide either with $\pbar'$ or with $\pbar''$.
To get the contradiction it is enough to show that
\begin{equation}\label{contr}
|p_0' - q| < 2R \qquad \hbox{and} \qquad |p_0''- q | < 2R;
\end{equation}
indeed, if (\ref{contr}) is true, $B(q, R)$ 
intersects both $B(p'_0, R)$ and $B(p''_0, R)$,
while it
must be contained either in $E$ 
or in the complement of $E$.
Let us compute the distance between $p_0'$ and $q$: 
\begin{eqnarray*}
|p_0'- q|^2 & = & 
(\xbar - R \sin \alpha(\pbar))^{2} + (\ybar + R \cos \alpha(\pbar) + R)^{2} \\ 
& = & 2 R^{2} + \xbar^{2} + \ybar^2 - 2R\, \xbar\sin \alpha(\pbar) 
+ 2R^2 \cos \alpha(\pbar) + 2R\, (1 + \cos \alpha(\pbar) )\ybar.
\end{eqnarray*}
Using the estimate $-R + \sqrt{R^2 - \xbar^2} \leq \ybar 
\leq R - \sqrt{R^2 - \xbar^2}$, the absurd assumption
and (\ref{sin}), we obtain
\begin{eqnarray*}
|p_0'- q|^2 & \leq & 
6 R^2 - 4R \sqrt{R^{2} - \xbar^2} - 2R\,\xbar \sin \alpha(\pbar) 
+2R (2R - \sqrt{R^{2} - \xbar^2}) \cos \alpha (\pbar)\\
& < & 6R^{2} - 4R \sqrt{R^{2} - \xbar^2} - 2\xbar^2 
+ 2 ( 2R -  \sqrt{R^2 - \xbar^2})\sqrt{R^2 - \xbar^2} \\
& = &  4R^2.
\end{eqnarray*}
By similar computations one can estimate the distance between
$p_0''$ and $q$.

Let $p=(\xbar, y)$ with $\ybar < y \leq R- \sqrt{R- \xbar^2}$
and let us suppose that
\begin{equation}\label{last}
\cos \alpha (\pbar) \geq \frac{\sqrt{R^2 - \xbar^2}}{R}. 
\end{equation}
We want to check that 
$$|p- q|^2 \leq R^{2},$$
which, by easy computations, is equivalent to
$$y - \ybar < 2 R \cos \alpha(\pbar).$$
By assumption we know that
\begin{eqnarray*}
y - \ybar & \leq & R - \sqrt{R^{2} - \xbar^{2}} - \ybar \\
& \leq & 2R - 2 \sqrt{R^{2} - \xbar^{2}} \\
& < & 2 \sqrt{R^{2} - \xbar^{2}}
\ \leq \ 2R \cos \alpha(\pbar),
\end{eqnarray*}
where the two last inequalities follow by the hypothesis 
$|\xbar | < \sqrt{3}R/2$ and by (\ref{last}).
\qed

Fixed $x_{1}$ in $[0, \sqrt{3}R/2[$, 
let us suppose by contradiction that the straight 
line $x=x_{1}$ intersects $\partial E 
\cap C^+$ in two distinct points $p_{1}$ and $q_{1}$. 
Then, if we call $p_{1}$ 
the point with smallest $y$-coordinate, 
by Lemma \ref{grafico} it follows that either $B(p_{1}', R)$ 
or $B(p_{1}'', R)$ must contain 
the point $q_{1}$ and this is impossible.
Therefore, we can conclude that $\partial E \cap C^+$ is the graph
of a function $f$. 

Since $f$ is between the functions $-R + \sqrt{R^{2} - x^{2}}$ and 
$R - \sqrt{R^{2} - x^{2}}$, 
which are both differentiable at $x=0$ with null derivative,
$f$ is differentiable at $x=0$ with derivative equal to $0$.
By a change of coordinates, we can repeat the same argument 
at every point belonging to $[0, \sqrt{3}R/2[$;
therefore,
$f$ is differentiable in $[0, \sqrt{3}R / 2[$ 
and the tangent line to the graph of $f$ at any point coincides with
the tangent line to the spheres associated to the same point.
From here, we obtain by Lemma \ref{grafico} the following bound on the 
norm of the derivative of $f$:
\begin{equation}\label{bound}
\left| \frac{df}{dx}(x) \right| \leq  \frac{|x|}{\sqrt{R^{2}- x^{2}}},
\end{equation}
for every $x \in [0, \sqrt{3}R/2[$.

To conclude the proof of the proposition in the case $n=2$, 
it is sufficient to check that 
the derivative of $f$ 
is Lipschitz with constant depending only on $R$.
First, we observe that, by (\ref{bound}),
\begin{equation}\label{bound2}
\left| \frac{df}{dx}(x) \right|  \leq 
\frac{2|x|}{R},
\end{equation}
for every $x\in [0, \sqrt{3}R/2 [$.
Given $p_{1}=(x_{1}, f(x_{1}))$ and $p_{2}=(x_{2}, f(x_{2}))$, 
we consider the following change of coordinates:
$$\begin{cases}
\tilde{x} \ =  \ \frac{1}{1+ \left( \frac{df}{dx}(x_{1})\right)^{2}} 
\left( x-x_{1}
+ (z - f(x_{1})) \frac{df}{dx}(x_{1}) \right) \\
\\
\tilde{z} \ = \ \frac{1}{1+ \left( \frac{df}{dx}
(x_{1})\right)^{2}} \left(- \frac{df}{dx} (x_{1}) (x-x_{1})
+ z - f(x_{1}) \right),
\end{cases}$$
which transform the point $p_{1}$ in the origin and the 
tangent line to $\partial E$ at $p_{1}$ in the $\tilde{x}$-axis.
With respect to the new coordinates, $\partial E$ 
is locally the graph of a function $\tilde{f}$ 
and the point $p_{2}$ has coordinates
$(\tilde{x}_{2}, \tilde{f}(\tilde{x}_{2}))$; then, by (\ref{bound2}),
\begin{equation}\label{boundtilde}
\left| \frac{d \tilde{f}}{d \tilde{x}} (\tilde{x}_{2})\right|
\leq \frac{2|\tilde{x}_{2}|}{R},
\end{equation}
if 
\begin{equation}\label{kappa}
|\tilde{x}_{2}| < \frac{\sqrt{3}}{2} R.
\end{equation} 
If we denote by $L$ 
the Lipschitz constant of $f$ in $[0, \sqrt{3}R/2 [$, 
we have that 
\begin{equation}\label{disttilde}
|\tilde{x}_{2}| \leq |p_{1} - p_{2}| \leq \sqrt{1 + L^{2}} |x_{1} - x_{2}|.
\end{equation}
Therefore, the condition (\ref{kappa}) is satisfied if 
$\displaystyle |x_{1} - x_{2}| \leq \frac{\sqrt{3} R}{2\sqrt{1+L^{2}}}
=:\lambda.$

By the relation 
$$\frac{df}{dx}(x_{2}) - \frac{df}{dx} (x_{1}) 
= \left(1 +  \frac{df}{dx}(x_{1})
\, \frac{df}{dx} (x_{2}) \right)\frac{d \tilde{f} }{d \tilde{x}} 
(\tilde{x}_{2}),$$
by (\ref{boundtilde}), and (\ref{disttilde}), it follows that
\begin{eqnarray}
\left| \frac{df}{dx}(x_{2}) - \frac{df}{dx} (x_{1}) \right| & \leq  & 
\left| 1 + 
\frac{df}{dx}(x_{1}) 
\frac{df}{dx}(x_{2}) \right| 
\left| \frac{d\tilde{f} }{d \tilde{x}}(\tilde{x}_{2}) 
\right| \nonumber \\
& \leq &  \frac{ 2\sqrt{1+L^2}}{R}\, \left| 1 + \frac{df}{dx}(x_{1})
\frac{df}{dx}(x_{2}) \right| |x_{1}-x_{2}|.
\end{eqnarray}
By the boundedness of the derivative of $f$, we can conclude that
there exists a positive constant
$c$, depending only on $R$, such that,
if $|x_{1} - x_{2}| \leq \lambda$,
then 
$$\| \frac{df}{dx}(x_{2}) - \frac{df}{dx}(x_{1}) \| \leq c |x_{1} - x_{2}|.$$

In the case
$|x_{1} - x_{2}| >\lambda$, we can find a finite number of points
$y_{0}:=x_{1} < y_{1} < \ldots < y_{k-1} < y_{k}:=x_{2}$ such that 
$|y_{j+1} - y_{j}| \leq \lambda$ for every $j=0, \ldots, k-1$.
Then, we obtain 
\begin{eqnarray*}
\left| \frac{df}{dx}(x_{1}) - \frac{df}{dx} (x_{2}) \right| & \leq & 
\sum_{j=0}^{k-1}
\left| \frac{df}{dx} (y_{j+1}) - \frac{df}{dx} (y_{j}) \right| \\
& \leq & c \sum_{j=0}^{k-1} |y_{j+1} - y_{j}| 
= c |x_{1} - x_{2}|.
\end{eqnarray*}
The proposition in the case $n=2$ is proved.

In the case $n\geq 3$ we can reduce to the $2$-dimensional one
by a slicing argument. For simplicity we sketch the proof
only for $n=3$; the general case can be treated in the same way.

From now on, we will write the coordinates of a point $p\in \rtre$
as a pair $(x,z)\in \rdue {\times } \reali$.
Given $E\in {\cal U}_R$ and 
$p\in \partial E$, we denote by $\proj{p}{E}$ the projection
on the plane 
which is tangent at $p$ to the balls $B(p', R)$ and $B(p'',R)$.

If $q\in \partial E$ we define
\begin{eqnarray*}
\avec{p}{E}(q)  & : = & \proj{p}{E} \left( \frac{1}{2R}(q''- q')  \right),\\
\bvec{p}{E}(q)  & :=  & \sqrt{1 - ( \avec{p}{E}(q))^2}.
\end{eqnarray*}

\begin{lemma}\label{primo}
There exist two constants $\delta>0$, $M>0$ such that,
for every $E\in {\cal U}_R$ and for every
$p,q \in \partial E$ with
$| \proj{p}{E} (p-q) | < \delta $, 
it results that
$\bvec{p}{E}(q) > M$. 
\end{lemma}  

{\sc Proof.} Let us suppose by contradiction that
for every $h\in \enne$ there exist $E_h \in {\cal U}_R$,
$p_h, q_h \in \partial E_h$ such that
\begin{equation}\label{ass}
| \proj{p_h}{E_h} (p_h - q_h )| \leq \frac{1}{h}, \qquad 
0\leq \bvec{p_h}{E_h}(q_h) \leq
\frac{1}{h}.
\end{equation}

Up to rototranslations, we can suppose
that $p_h =(0,0)$, $p'_h=(0, -R)$, and $p''_h =(0,R)$.
If we denote by $(x_h, z_h)$ the coordinates of $q_h$, we obtain
that
$$| p_h'' - q_h'|^2 = |(x_h,0) + R\, \avec{p_h}{E_h}(q_h)|^2 +
(z_h - R\, b_h(q_h) -R)^2.$$
Since by (\ref{ass}) the right-hand side tends to $2R^2$ as $h\to \infty$,
for $h$ large the ball $B(p''_h, R)$ intersects $B(q_h', R)$,
which is impossible. \qed

Now we are in position to prove the crucial lemma which allows us to perform the two-dimensional reduction.

\begin{lemma}\label{secondo}
Let $\delta>0$ and $M>0$ as in Lemma \ref{primo}.
Let $E$ be in ${\cal U}_R$,  $p\in \partial E$ and choose 
a coordinate system in $\esse{p}{E}$. Then,
for every $(\xbar,0)$ with $|\xbar |< \delta$
the section of $E$ with any vertical plane $\gamma$ passing through
$(\xbar,0)$ satisfies in $\gamma$ the exterior and interior sphere 
condition
with radius $MR$ at every point of $\partial E \cap C$,
where $C:=
\{ x\in\reali^{n-1}: |x|< \delta  \} { \times } ]{ - } R, R[$. 
\end{lemma}

{\sc Proof.} Let $\gamma$ be a vertical plane passing through
$(\xbar,0)$ and let $(v,0)$ be a unit normal vector to $\gamma$.
Let $q\in \partial E \cap C \cap \gamma$.
By Lemma \ref{primo}, we have that
$$|a_p^E(q)| = \sqrt{1 - (b_p^E(q))^2} < \sqrt{1-M^2};$$
hence, if we call $\alpha$ the angle in $[0, \pi[$ between
$q''-q'$ and $(v,0)$, then
\begin{equation}\label{alfa}
|\cos \alpha| =|a_p^E(q)\cdot (v,0)|< \sqrt{1-M^2}.
\end{equation}
Then the point $q$ satisfies the exterior and interior 
sphere condition
in
$\gamma$ with radius $R\sin \alpha$, which by (\ref{alfa})
is greater than $MR$. 
\qed

Now we can prove the proposition in the case $n=3$.

Given $E\in {\cal U}_R$ and $p\in \partial E$, we choose
a coordinate
system in $\esse{p}{E}$ and we call $C$ the cylinder
$\{ x\in \rdue: |x| < \overline{\delta} \} { \times } ] {-}R, R[$,
where $\overline{\delta} := \min \{ \delta, \sqrt{3}MR/2 \}$, and
$\delta$, $M$ are as in Lemma \ref{primo}. 
Applying the $2$-dimensional result to the
sections of $E$ with the vertical planes passing through
the point $p$, by Lemma \ref{secondo} we obtain that
$\partial E \cap C$ is the graph of a function $f$ defined in $\{ x\in \rdue: |x| < \overline{\delta} \}$ .

To show the differentiability of $f$, we can repeat the same argument as in
the $2$-dimensional case. Moreover, Lemma \ref{primo} gives an
uniform bound on
the norm of the gradient of $f$.

Using Lemma \ref{secondo}, the 
$2$-dimensional result, and Lemma \ref{primo}, 
we can find $\rho  \in ]0, \overline{\delta}]$
and $N>0$ such that in $\{ x\in \rdue: |x| < \rho \}$
the restriction of $f$ to any straight line is a function of class 
$\sobdue$ with $\sobdue$-norm less than $N$.

To conclude, we define the function 
$$g(x_1, x_2):= \lim_{h\to \infty} h \left[ \partial_{x_1} f\left( x_1+\frac{1}{h},x_2 \right) 
-\partial_{x_1} f(x_1,x_2) \right],$$
for a.e. $x=(x_1,x_2)$.
By the above remark, $g$ is defined a.e. and belongs to $L^{\infty}$ with
$L^{\infty}$-norm less than $N$. 
Using the absolute continuity of $\partial_{x_1} f$ on
the straight lines $x_2={\rm constant}$, it is easy to check that
$g$ coincides with the second distributional derivative $\partial^2_{x_1} f$ .
Analogously, we can prove that there exists $\partial^2_{x_2} f$ in the
distributional sense, and that it belongs to $L^{\infty}$ with
$L^{\infty}$-norm less than $N$. 
To show that $\partial_{x_1} \partial_{x_2} f$ exists and belongs to
$L^{\infty}$ with
$L^{\infty}$-norm less than $N$, one can argue in a similar way, 
by considering
the restriction of $f$ to the straight lines $x_1 - x_2=
{\rm constant}$.
\qed

\begin{lemma}\label{diametri}
Let $\{ E_h \}_h$ be a sequence of connected sets in ${\cal U}_R$
such that $\lim_{h\to \infty}\diam (E_h) =+ \infty$. Then
$$\lim_{h\to \infty} {\cal L}^n (E_h) = +\infty.$$
\end{lemma}

{\sc Proof.}
Since $\lim_{h\to \infty}\diam (E_h) =+ \infty$, 
for every $h\in \enne$ we can find
$p_1^h, \ldots, p_{m_h}^h \in \partial E_h$,
where $m_h$ is the integer part of 
$\diam(E_h) / 4R$, such that
$| p_i^h - p_j^h| \geq 4R$ for every $i\neq j$.
We clearly 
have that $\{B((p_{h}^{i})',R)\}_{i=1,\ldots, m_h}$  
is a family of disjoint balls all contained in $E_h$; hence,
$${\cal L}^{n}(E_h)\geq m_h {\cal L}^{n}(B(0,R)),$$
and the second term goes to infinity as
$h\to \infty$. \qed

\section{The compactness result}

In the sequel, if $\{ f_j \}_j$ is a sequence in
$\sobdue(\Omega)$ and $f$ is a function in $\sobdue(\Omega)$, we mean
by the notation $f_j  \rightharpoonup f$ in w$^*$-$\sobdue(\Omega)$
that the sequence $\{ f_j\}_j$ converge to $f$ in the weak$^*$- topology
of $\sobdue(\Omega)$.
Given $E\subset \rn$, we denote the characteristic 
function of $E$ by $\chi_E$. If $\partial E$ is sufficiently regular,
we denote the unit outer normal vector to $\partial E$ at the point $p$
by $\nu_{\partial E}(p)$.

We start by  recalling  two notions of set-convergence.

\begin{definition}
Let $\{E_h \}_h$ and $E$ be measurable 
subsets of $\reali^{n}$. We say 
that the sequence $\{ E_{h} \}_{h}$ converges to $E$ a.e.
if $\chi_{E_h} \to \chi_E$ a.e., and that
$\{ E_{h} \}_{h}$ converges to $E$
in $L^1$ if 
$\chi_{E_h} \to \chi_E$ in $L^1(\rn )$.
\end{definition}

\begin{definition}\label{curat}
Let $\{ E_{h} \}_h$ and $E$ be closed subsets of $\reali^{n}$. We say 
that the sequence $\{ E_{h} \}_{h}$ converges to $E$ 
in the sense of Kuratowski (and we write  
$E_{h} \kconv E$) if 
\begin{description}
\item{i)} $p_{h}\in E_{h},\ \exists\, p_{h_{k}}\to p\ \Rightarrow \ p\in 
E$;
\item{ii)} $\forall p\in E,\ \exists\, p_{h}\in E_{h}\ :\ p_{h}\to p.$
\end{description}
\end{definition}
It is well known that on the space of equibounded compact sets, the
Kuratowski 
convergence is induced by the Hausdorff distance.

\begin{theorem}\label{cpt}
Let $\{ E_h  \}_h$ be an equibounded sequence of sets 
belonging to ${\cal U}_R$.
Then there exist $E\in {\cal U}_R$ and a subsequence 
$\{ E_{h_{j}} \}_j$ such that
\begin{description}
\item[a)] $E_{h_{j}} \kconv E$ and
	$E_{h_{j}} \to E$ in $L^1$;
\item[b)] $\partial E_{h_{j}} \kconv \partial E$ and
	$\lim_{j\to \infty} 
	{\cal H}^{n-1}( \partial E_{h_{j}})= {\cal H}^{n-1}
	(\partial E)$;
\item[c)] there exists a constant $\eta \in ]0,1[$
(depending only on $R$), such that for every $p\in \partial E$,
if we call $C^{\eta}$ the cylinder
$\{ x\in \reali^{n-1}: |x| \leq \eta R \} 
{\times} [{-}\eta R, \eta R]$ expressed with 
respect to any  coordinate system belonging to $\esse{p}{E}$,
 and $S^{\eta}$ the section $C^{\eta} \cap \{ z=0 \}$,
then $\partial E \cap C^{\eta}$ is the graph 
of a function $f\in \sobdue(S^{\eta})$, and
$\partial E_{h_{j}} \cap C^{\eta}$ is definitively
the graph of a function $f_{j} \in \sobdue(S^{\eta})$.
Moreover, $f_j \rightharpoonup f$ 
in w$^{*}$-$\sobdue (S^{\eta} )$.
\end{description}
\end{theorem}

{\sc Proof.} 
Since $\{ E_h \}_h$ is equibounded,
there exist a compact set $E$ and 
a subsequence, which we denote again by $\{ E_{h} \}_{h}$, such that
\begin{equation} \label{estrazione}
E_{h} \kconv E.
\end{equation}
Let us prove that $E\in {\cal U}_R$.

First of all, we remark that 
if $\{ p_h \}_h$ is a sequence such that $\dist ( p_{h}, E_{h}) > c > 0$ for every $h\in \enne$, then
every limit point $p$ of $\{ p_h \}_h$ belongs to $\complement E$.
Indeed, let us suppose by contradiction that
there exists $\{ p_{h_{k}} \}_k$ which converges to $p \in E$; then, by ii) in 
Definition \ref{curat}, for every $h \in \enne$ there is 
$q_{h}\in E_{h}$ 
such that $\{ q_{h}\}_h$ converges to $p$ and so, 
$| q_{h_{k}} - p_{h_{k}}| \to 0$, in contradiction with the initial assumption.

{\sc Claim 1.} Every point $p\in \partial E$ is the limit of a sequence $\{ p_h \}_h$ such that
$p_h \in \partial E_h$ for every $h\in \enne$.

Let $p\in \partial E$. By ii) in Definition \ref{curat} there 
exists $p_{h}\in E_{h}$ such that
$\{ p_{h} \}_h$
converges to $p$; clearly, it is 
enough to show that
$\dist (p_{h},\partial E_{h}) \to 0$.
If by contradiction 
there exists a subsequence $\{ p_{h_{k}}\}_{k}$ 
such that $\dist(p_{h_{k}},\partial 
E_{h_{k}}) > c>0$ for every $k\in \enne$, 
then the ball $B(p_{h_{k}},c)$ is contained in
$E_{h_{k}}$. Since for every $q\in B(p,c)$ we can find 
$q_k \in B(p_{h_{k}}, c)$ such that 
$q_k \to q$, then by i) in Definition \ref{curat},
$q\in E$. Therefore $B(p,c)\subset E$, hence $p \in  
\inter{E}$, which contradicts our initial 
assumption.

{\sc Claim 2.} If $p_{h}\in \partial E_{h}$ for every $h\in \enne$ and 
there is a subsequence 
$\{ p_{h_{k}} \}_k$ converging to a point $p$,
then $p\in \partial E$.
Moreover, there exist $p'$, $p''$ such that
$$p \in  \partial B(p', R) \cap \partial B(p'', R), 
\quad  B(p', R) \subset E, \quad  
B(p'', R) \cap E = \emptyset.$$

Since $E_{h_{k}} \in {\cal U}_R$ for every $k\in \enne$,
there exist $p'_k$, $p''_k$ such that
the balls $B(p_{k}',R)$, $B(p''_{k},R)$ are
contained respectively in $E_{h_{k}}$ and in
$\complement E_{h_{k}}$.
Up to subsequences, we can suppose that $\{ p'_{k} \}_k$
and $\{ p''_{k} \}_k$
converge to $p'$ and $p''$
respectively.
Therefore,
$$\overline{B(p'_{k}, R)} \, \kconv \, \overline{B(p',R)}, \qquad  
\overline{B(p''_{k}, R)} \, \kconv
\, \overline{B(p'',R)},$$
and, since $\{ p_{h_{k}} \} =  \partial B(p'_k,R) \cap \partial B(p''_k,R)$,
we have that
\begin{equation} \label{intersezione}
\{ p \} =  \partial B(p',R) \cap \partial B(p'',R).
\end{equation}

If $q \in B(p',R)$, then $q$ is the limit of a sequence 
$\{ q_{k}\}_k$ such that $q_k \in B(p'_{k}, R)\subset E_{h_{k}}$; by (\ref{estrazione}) 
and i) in Definition \ref{curat}, it follows that $q\in E$; this means that 
$B(p',R)$ is contained in $E$.

Let $q\in B(p'',R)$ and let $q_{k}:= q-p'' -p''_{k}$
for every $k\in \enne$.
It is clear that
$q_{k} \in B(p''_{k},R)$, 
there exists a constant $c>0$ such 
that $\dist (q_{k},E_{h_{k}}) = c$, and  the
sequence $\{q_{k}\}_k$ converges to $q$. Thus, as
remarked before,
$q \in \complement E$. We can conclude that
$B(p'',R)$ is contained in the complement of $E$. 

By (\ref{intersezione}), it follows that $p\in \partial E$ and this concludes the proof
of the claim.

By Claim 1 and 2, we can deduce that $E\in {\cal U}_R$ and also 
\begin{equation} \label{convbord}
\partial E_{h} \kconv \partial E.
\end{equation}

To show the convergence in $L^1$, 
it is enough to prove the pointwise convergence of
$\{ \chi_{E_{h}} \}_h$ to $\chi_E$ for every $p\notin \partial E$; 
indeed, by the regularity of
$E$, we have that ${\cal L}^n (\partial E)=0$.
If $p\in \inter{E}$, then by (\ref{estrazione}) and (\ref{convbord}) 
there exists $p_h \in \inter{E_h}$ such that 
$\dist (p_h, \partial E_h) > c>0$ and $p_h \to p$. Then $p$ 
definitively belongs to $B(p_h, c)$, 
which is contained in $\inter{E_h}$; hence,
$\chi_{E_{h}}(p) =1$ for $h$ large and so,
$\{ \chi_{E_{h}}(p)\}_h$ obviously converges to $\chi_E (p)$.
If $p\in \complement E$ and by contradiction 
there exists a subsequence $\{ h_k \}_k$
such that $p \in E_{h_{k}}$, then by i) in Definition \ref{curat} 
$p\in E$, which is absurd.

Let us prove the third part of the proposition.

Let $p\in \partial E$. By (\ref{convbord}), there is a sequence $\{ p_h \}_h$
such that $p_h \in \partial E_h$ for every $h\in \enne$ and $p_h \to p$.
From now on, we will work in a coordinate system
belonging to $\esse{p}{E}$.
By Proposition \ref{regular}, there exists $\delta \in ]0, 1[$,
depending only on $R$, 
such that, if we set
$C:= \{ x\in \reali^{n-1}: |x|< \delta R \} { \times }
]{ - } R , R [$, then $\partial E \cap C$ 
is the graph of a function $f$ defined on the base
of $C$ and of class $\sobdue$. 
Let us denote by $C_h$ the cylinder obtained by translating 
the centre of $C$ in $p_h$
and by rotating the axis of $C$ 
in such a way that it is directed along
$\nu_{\partial E_h}(p_h)$.
By Proposition \ref{regular}, 
$\partial E_h \cap C_h$ is the graph of a function $f_h$ defined 
on the base of $C_h$ and of class $\sobdue$.
We recall that
\begin{equation}\label{equi}
\| f_h \|_{\infty} \leq (1 - \sqrt{1-\delta^2})R,
\end{equation}
(see the proof of Proposition \ref{regular}).
Since $\nu_{\partial E_h}(p_h)$ is parallel to the vector $p''_h - p'_h$, 
$\{ p'_h \}_h$ converges
to $p'$, $\{ p''_h \}_h$ converges to $p''$ (see the proof of Claim 2),
and $\nu_{\partial E}(p)$ is parallel to the vector $p'' - p'$,
we have that
\begin{equation}\label{normali}
\nu_{\partial E_h}(p_h) \to \nu_{\partial E}(p).
\end{equation}
By the convergence of $\{ p_h \}_h$ to $p$ and by (\ref{normali}), 
it follows that 
for $h$ sufficiently large $C_h$ contains the cylinder $C^{\eta}
= \{ x\in \reali^{n-1}: |x| \leq \eta R \} {\times} [{-}\eta R, \eta R]$, where
$\eta \in ] 1-\sqrt{1-\delta^2}, \delta[$.
Using (\ref{normali}), (\ref{equi}) and the equiboundedness of
$\{ \nabla f_h \}_h$, one can easily check that for $h$ large enough
$\partial E_h \cap C^{\eta}$ can be expressed as the graph of a new function 
$\tilde{f}_h$
defined on the base of $C^{\eta}$.

Using again (\ref{normali}) and the equiboundedness of $\{ f_h \}_h$ in
$\sobdue$-norm, it is easy to see that 
$\tilde{f}_h \in \sobdue(S^{\eta})$
and the $\sobdue$-norm of $\tilde{f}_h$ is bounded 
by a constant depending only on $R$.
Then there exist a subsequence $\{\tilde{f}_{h_{k}} \}_k$ and a function
$\tilde{f} \in \sobdue(S^{\eta })$ such that
$\{ \tilde{f}_{h_{k}} \}_k$ converge to $\tilde{f}$ 
in w$^{*}$-$\sobdue(S^{\eta})$
(and then in $C^{1}$-norm).
It remains to prove that $\tilde{f}$ coincides with $f$ on $S^{\eta}$.

{\sc Claim 3.} It results that
\begin{equation}
\graph \tilde{f}_{h_{k}} \kconv \graph \tilde{f}
\end{equation}
and
\begin{equation}
\partial E_h \cap C^{\eta} \kconv \partial E \cap C^{\eta}.
\end{equation}

Let $p_{k} \in \graph \tilde{f}_{h_{k}}$ and let $\{ p_{k_{j}} \}_j$ 
be a subsequence converging to a point
$p$. The point $p_k$ has coordinates $(x_k, \tilde{f}_{h_{k}}(x_k))$ with 
$|x_k|\leq \eta R$; up to subsequences, $\{x_{k}\}_k$ converges to 
a point $x$ such that $|x|\leq \eta R$. By the uniform
convergence of the functions, we obtain that
$\{ p_k\}_k$ tends to the point $(x, \tilde{f}(x))$, which belongs
trivially to $\graph \tilde{f}$.
Then property i) in
Definition \ref{curat} is proved. 
Let $p=(x, \tilde{f}(x)) \in \graph \tilde{f}$ with $|x|\leq \eta R$.
The point $p_k := (x, \tilde{f}_{h_{k}}(x))$ belongs 
to $\graph \tilde{f}_{h_{k}}$ and $\{ p_k \}_k$ converges
to $p$; hence, property ii) in Definition \ref{curat} is verified.

Since $C^{\eta}$ is closed and by (\ref{convbord}), 
property i) in Definition \ref{curat} 
is trivial.
By (\ref{convbord}) property ii) is easily verified for the points 
belonging to $\partial E \cap \inter{C^{\eta}}$;
if $p\in \partial E \cap \partial C^{\eta}$ and $p=(x,z)$ with 
$|x|=\eta R$ and $|z|\leq \eta R$,
then it is enough to take the sequence $p_h= (x, \tilde{f}_h(x))\in \partial
E_h \cap \partial C^{\eta}$. 

By Claim 3,
since $\graph \tilde{f}_{h_{k}} = \partial E_{h_{k}} \cap C^{\eta}$,
it follows that $\graph \tilde{f}$ coincides with $\partial E\cap C^{\eta}$.
Then, $\tilde{f} =f$ on $C^{\eta}$ and the whole sequence 
$\{ \tilde{f}_h \}_h$ converges
to $f$ in w$^*$-$\sobdue(S^{\eta})$.

Let us prove the second part of b).

By point c), for every $p\in \partial E$
there exists a cylinder $C$ centred at $p$, with base a 
$(n-1)$-dimensional sphere $S$, such that
$\partial E \cap C$ is the graph of a function $f\in \sobdue(S)$,
for $h$ large
$\partial E_h \cap C$ is the graph of a function $f_h\in \sobdue(S)$,
and $f_h \rightharpoonup f$ in w$^*$-$\sobdue(S)$.
We can recover $\partial E$ with a finite number of these cylinders
$C_1, \ldots, C_m$. Let us call $f^i_h$ the function such that
$\graph f^i_h = \partial E_h \cap C_i$, and $f^i$ the function such that
$\graph f^i = \partial E \cap C_i$.

Let $\varepsilon >0$ be such that 
$$(\partial E)_{\varepsilon} := \{ p\in \rn : 
\dist( p, \partial E) \leq \varepsilon \}
\subset \cup_{i=1}^m C_i.$$
We can consider a partition of unity associated to the
recovering $\{ C_1, \ldots, C_m \}$,
i.e. a family of functions 
$\phi_i \in C^{\infty}_0(C_i)$ ($i=1, \ldots, m$)
such that $0 \leq \phi_i \leq 1$,
$$\sum_{i=1}^{m} \phi_i =1 \ \hbox{on} \ (\partial E)_{\varepsilon}, \qquad
\sum_{i=1}^{m} \phi_i \leq 1 \ \hbox{on} \ \cup_{i=1}^m C_i.$$
By (\ref{convbord}), for $h$ large $\partial E_h \subset
(\partial E)_{\varepsilon}$. Then,
$${\cal H}^{n-1}(\partial E_h ) =
\sum_{i=1}^m \int_{\partial E_h\cap C_i} \phi_i\, d{\cal H}^{n-1} =
\sum_{i=1}^m \int_{\graph f^i_h} \phi_i\,
 d{\cal H}^{n-1}.$$
Using the Area Formula and the $C^1$-convergence of
$\{ f^i_h \}_h$ to $f^i$, it is easy to see that
for every $i=1, \ldots, m$
$$\lim_{h\to \infty} \int_{\graph f^i_h} \phi_i\,
 d{\cal H}^{n-1} 
= \int_{\graph f^i} \phi_i\,
d{\cal H}^{n-1}.$$
Therefore,
$${\cal H}^{n-1}(\partial E )  =   \sum_{i=1}^m \int_{\graph f^i} \phi_i\,
d{\cal H}^{n-1} 
=   \lim_{h\to \infty} \sum_{i=1}^m \int_{\graph f^i_h} \phi_i\,
d{\cal H}^{n-1} 
 =  \lim_{h\to \infty} {\cal H}^{n-1}( \partial E_h).$$
\qed

\section{The semicontinuity result}

Given $E\in {\cal U}_R$, we think $\partial E$
oriented by the outer normal field 
(all the results we will state still remain true
if we choose the opposite orientation). 
We denote the principal curvatures 
(i.e. the eigenvalues of the second fundamental 
quadratic form) of $\partial E$ at the point $x$
by $\kappa_i(x)$ with $i=1, \ldots, n-1$, and the $p^{th}$-elementary 
symmetric function of the principal
curvatures, called $p^{th}$-elementary symmetric curvature, by 
\begin{equation}\label{symm}
K_p(x) = \binom{n-1}{p}^{-1} \sum_{1\leq i_1 < \ldots < i_p \leq n-1} 
\kappa_{i_{1}}(x) \ldots \kappa_{i_{p}}(x)
\end{equation}
for $p= 1, \ldots, n-1$. We also use the notation
$$H:=K_1 \qquad \hbox{and} \qquad K:=K_{n-1}$$
for the mean curvature and the Gauss curvature respectively.
In the case $n=2$ we simply denote the curvature by $\kappa$.

It is well known from differential geometry (see \cite{spi})
that the $p^{th}$-elementary symmetric curvature
is the coefficient of the term of degree $n-1-p$
of the characteristic polynomial of the second
fundamental quadratic form.
If $\partial E$ is locally the graph of a function $f$,
then the second fundamental quadratic form is given by
the product $G^{-1}B$, where $G=(g_{ij})$ is the matrix
defined by
$$g_{ij} = \begin{cases} 
1 + (\partial_{x_i} f)^2 & \text{if $i=j$}, \\
\partial{x_i} f \, \partial{x_j} f & \text{if $i\neq j$},
\end{cases}$$
while $B=(b_{ij})$ is the matrix
$$b_{ij} = \frac{\partial_{x_i x_j} f}{\sqrt{1 + |\nabla f|^2}}.$$
By induction, it is easy to prove that for every $p=1, \ldots, n-1$
there exists a continuous function $\psi_p= \psi_p(s, \zeta)$,
linear with respect to $\zeta$, such that
$$K_p(x, f(x))= \psi_p(\nabla f(x), M(\nabla^2 f(x))),$$ 
where $M(\nabla^2 f(x))$ is the vector of the determinants of all the minors
of $\nabla^2 f(x)$.

In the sequel we will consider functionals of the form
$$F(E):= \int_{\partial E} \varphi( K_1, \ldots, K_{n-1}) \, d{\cal H}^{n-1},$$
where $\varphi: \reali^{n-1}\to \reali$ is a 
given convex function.
Functionals of this type arise in different contexts; for instance:
\begin{itemize}
\item the Willmore's functional (see \cite{ch, wil1}),
$\displaystyle F(E)= \int_{\partial E} |H(x)|^{n-1} d{\cal H}^{n-1}(x);$
\item $\displaystyle F(E) = \int_{\partial E} [H^2 - K_2]^{(n-1)/2} 
d{\cal H}^{n-1},$
studied in \cite{wil2};
\item $\displaystyle
F(E)= \int_{\partial E} \varphi
\left( \sum_{i=1}^{n-1}\kappa_i^2(x) \right) 
d{\cal H}^{n-1}(x);$
in the case $n=2$ and $\varphi(x)=1+x$, we find the functional
considered in \cite{bel-dm-pao}:
$$F(E)=\int_{\partial E}(1 + \kappa^2)\, d{\cal H}^1.$$
\end{itemize}

\begin{theorem}\label{lsc}
Let $\varphi: \reali^{n-1} \to \reali$ be a convex function.
If $E\in {\cal U}_R$ and $\{ E_h \}_h$ is a sequence 
in ${\cal U}_R$ such that
$E_h \to E$ in $L^1$,
then 
$$\int_{\partial E} \varphi( K_1, \ldots, K_{n-1}) \, d{\cal H}^{n-1}
\leq \liminf_{h\to \infty} \int_{\partial E_h} \varphi(K_1, \ldots, K_{n-1}) 
\, d{\cal H}^{n-1}.$$
\end{theorem}

For the proof of the theorem we need the following lemma.

\begin{lemma}\label{sketch}
Let $\phi: \rn {\times }\reali^{n-1} \to [0, +\infty[$ be 
globally continuous and convex in the last $n-1$ variables.
Let $\Omega$ be an open bounded
subset of $\reali^{n-1}$ and let $f_h, f\in \sobdue(\Omega)$.
If $f_h \rightharpoonup f$ in w$^*$-$\sobdue(\Omega)$, then
$$\int_{\graph f} \phi(q, K_1, \ldots, K_{n-1}) \, d{\cal H}^{n-1}(q)
\leq \liminf_{h\to \infty} \int_{\graph f_h} \varphi(q,K_1, \ldots, K_{n-1}) 
\, d{\cal H}^{n-1}(q).$$
\end{lemma}

{\sc Proof.}
It is not restrictive to assume that $\Omega$ is smooth.

As remarked above,
for every $p=1, \ldots, n-1$ and for every $x\in \Omega$ we have that
$$K_p (x, f_h(x)) = \psi_p (\nabla f_h(x), M(\nabla^2 f_h(x))),$$
where $\psi_p$ is globally continuous and linear in the second variable.

Using the Area Formula, we can write
\begin{equation}\label{loc}
\int_{\graph f_h} \phi(q, K_1, \ldots, K_{n-1}) 
\, d{\cal H}^{n-1}(q) = \int_{\Omega} \phi'(x, f_h(x),\nabla f_h(x),
\nabla^2 f_h(x)) \, dx,
\end{equation}
where 
$$\phi'(x,z,s,\xi) := \phi
((x,z), \psi_1(s, M(\xi)), 
\ldots, \psi_{n-1}( s, M(\xi)) ) \sqrt{1 + |s|^2}$$
for every $x\in \Omega$, $z\in \reali$, $s\in \rn$, and $\xi \in \matr$.
Let us define the function
$$\phi''(x,s,\xi) := \phi'(x,f(x), s,\xi)$$
for every $x\in \Omega$, $s\in \rn$, and $\xi \in \matr$.
Since $\phi''$ is positive, globally
continuous and polyconvex in $\xi$,
by Theorem II.1 in \cite{ace-fus}, it follows that
$$\int_{\Omega} \phi''(x, \nabla f(x),
\nabla^2 f(x)) \, dx \leq 
 \liminf_{h\to \infty} \int_{\Omega} \phi''(x, \nabla f_h(x),
\nabla^2 f_h(x)) \, dx,$$
that is
\begin{equation}\label{liminf}
\int_{\Omega} \phi'(x, f(x), \nabla f (x),
\nabla^2 f (x)) \, dx \leq 
\liminf_{h\to \infty} \int_{\Omega} \phi'(x, f(x), \nabla f_h (x),
\nabla^2 f_h (x)) \, dx .
\end{equation}
Using the uniform continuity of $\phi'$ on bounded sets
and the uniform convergence of $\{ f_h \}_h$ to $f$,
we have that
\begin{equation}\label{bb}
\lim_{h \to \infty} \int_{\Omega} 
|\phi' (x, f_h(x), \nabla f_h(x),
\nabla^2 f_h(x)) - \phi' (x, f(x), \nabla f_h(x),
\nabla^2 f_h(x))| \, dx =0.
\end{equation}
By (\ref{loc}), (\ref{liminf}), and (\ref{bb}),
the thesis easily follows.
\qed

{\sc Proof of Theorem \ref{lsc}.}
First of all, we observe that the sequence $\{ E_h \}_h$ is equibounded;
indeed, let $M>0$ be such that $E\subset B(0,M)$ and let
$\tilde{E}_h$ be the union of all the connected components
of $E_h$ which intersect $B(0,M)$.
By the $L^1$-convergence of $\{ E_h \}_h$ to $E$, it is clear
that
$$\lim_{h\to \infty} {\cal L}^n (E_h \setminus \tilde{E}_h) =0.$$
Recalling that, if $E\neq \emptyset$ belongs to ${\cal U}_R$,
then ${\cal L}^n (E) \geq {\cal L}^n(B(0,R))$,
we deduce that for $h$ large $E_h = \tilde{E}_h$, i.e.
all the connected components of $E_h$ definitively intersect 
$B(0,M)$.
The equiboundedness easily follows by Lemma \ref{diametri}, and
allows to conclude that the sequence $\{ E_h \}_h$
satisfies a), b), c) of Theorem \ref{cpt}.
(Note that we have incidentally proved that
in the class ${\cal U}_R$, $L^1$-convergence
and Kuratowski convergence are actually equivalent).

Let us suppose for the moment that $\varphi$ is positive.
By Theorem \ref{cpt}, for every $p\in \partial E$
there exists a cylinder $C$ centred at $p$, with base a 
$(n-1)$-dimensional sphere $S$, such that
$\partial E \cap C$ is the graph of a function $f\in \sobdue(S)$,
for $h$ large
$\partial E_h \cap C$ is the graph of a function $f_h\in \sobdue(S)$,
and $f_h \rightharpoonup f$ in w$^*$-$\sobdue(S)$.
We can recover $\partial E$ with a finite number of these cylinders
$C_1, \ldots, C_m$. Let us call $f^i_h$ the function such that
$\graph f^i_h = \partial E_h \cap C_i$, and $f^i$ the function such that
$\graph f^i = \partial E \cap C_i$.

We can consider a partition of unity associated to the
recovering $\{ C_1, \ldots, C_m \}$,
i.e. a family of functions 
$\phi_i \in C^{\infty}_0(C_i)$ ($i=1, \ldots, m$)
such that $0 \leq \phi_i \leq 1$,
\begin{equation}\label{unity}
\sum_{i=1}^{m} \phi_i =1 \ \hbox{on} \ \partial E, \qquad
\sum_{i=1}^{m} \phi_i \leq 1 \ \hbox{on} \ \cup_{i=1}^m C_i.
\end{equation}
Then
\begin{eqnarray*}
\int_{\partial E} \varphi( K_1, \ldots, K_{n-1}) \, d{\cal H}^{n-1} & = &
\sum_{i=1}^m \int_{\partial E\cap C_i} \phi_i\,
\varphi( K_1, \ldots, K_{n-1}) \, d{\cal H}^{n-1} \\
& = & 
\sum_{i=1}^m \int_{\graph f^i} \phi_i\,
\varphi( K_1, \ldots, K_{n-1}) \, d{\cal H}^{n-1} \\
& \leq & \liminf_{h\to \infty} \sum_{i=1}^m \int_{\graph f^i_h} \phi_i\,
\varphi( K_1, \ldots, K_{n-1}) \, d{\cal H}^{n-1} \\
& = & \liminf_{h\to \infty} \sum_{i=1}^m \int_{\partial E_h \cap C_i} \phi_i\,
\varphi( K_1, \ldots, K_{n-1}) \, d{\cal H}^{n-1} \\
& \leq & \liminf_{h\to \infty} \int_{\partial E_h}
\varphi( K_1, \ldots, K_{n-1}) \, d{\cal H}^{n-1},
\end{eqnarray*}
where we used Lemma \ref{sketch} and (\ref{unity}). 

If $\varphi$ is bounded from below by a constant $c\in \reali$,
we can apply the previous argument to the function
$\varphi - c$, to conclude that
\begin{multline}\nonumber
\int_{\partial E} \varphi( K_1, \ldots, K_{n-1}) \, d{\cal H}^{n-1} 
- c {\cal H}^{n-1}(\partial E)    \leq  \\
 \leq  \liminf_{h\to \infty} \left( \int_{\partial E_h}
\varphi( K_1, \ldots, K_{n-1}) \, d{\cal H}^{n-1}
- c {\cal H}^{n-1}(\partial E_h) \right) \\
=  \liminf_{h\to \infty}  \int_{\partial E_h}
\varphi( K_1, \ldots, K_{n-1}) \, d{\cal H}^{n-1}
-c {\cal H}^{n-1}(\partial E),
\end{multline}
where we used property b) in Theorem \ref{cpt}.

Finally, if $\varphi$ is a generic convex function, let us set
$$c:= \inf \left\{  \varphi( K_1(p), \ldots, K_{n-1}(p)): 
p\in  \bigcup_{h\in \enne} \partial E_h \cup\partial E \right\},$$
which is finite by the equiboundedness of curvatures
(see Proposition \ref{regular}).
If we define $\tilde{\varphi} := \varphi \vee c$,
we have that $\tilde{\varphi}$ is a convex function bounded from below; hence,
\begin{eqnarray*}
\int_{\partial E}
\varphi( K_1, \ldots, K_{n-1}) \, d{\cal H}^{n-1}
& = & \int_{\partial E}
\tilde{ \varphi}( K_1, \ldots, K_{n-1}) \, d{\cal H}^{n-1}\\
& \leq & \liminf_{h\to \infty}  \int_{\partial E_h}
\tilde{\varphi} ( K_1, \ldots, K_{n-1}) \, d{\cal H}^{n-1}\\
& = & \liminf_{h\to \infty}  \int_{\partial E_h}
\varphi( K_1, \ldots, K_{n-1}) \, d{\cal H}^{n-1}.
\end{eqnarray*}
\qed

In the following proposition, we study the
asymptotic behaviour of $F$ when $R$ goes to $0$, in
the case $n=2$ and $\varphi(\kappa )=
1 + | \kappa |^p$, showing the relationship 
with the relaxed functional introduced in \cite{bel-dm-pao}.

\begin{proposition}\label{gamma}

Let $\{ F_{R} \}_{R>0}$ be the family of functionals
$$F_{R}(E) := \begin{cases} 
\int_{\partial E} (1 + |\kappa |^{p})\, d{\cal H}^{1} 
& \text{if $E \in {\cal M} \cap {\cal U}_{R}$}, \\
+ \infty & \text{otherwise in ${\cal M}$},  
\end{cases}$$
where $p> 1$, and ${\cal M}$ is the class of measurable bounded sets in $\rdue$.
Then, $F_{R}$, as $R\to 0$, $\Gamma$-converges 
(for the definition and the properties of $\Gamma$-convergence,
see \cite{dm})
with respect to the $L^{1}$-topology
to the lower semicontinuous envelope, $\overline{F_0}$, of 
$$F_0(E): = \begin{cases} 
\int_{\partial E} (1 + |\kappa|^{p})\, d{\cal H}^{1} 
& \text{if $E \in{\cal M}\cap  C^{2}$}, \\
+ \infty & \text{otherwise in ${\cal M}$}. 
\end{cases}$$
\end{proposition}

{\sc Proof.}
For every $E \in {\cal M}$, $\{ E_{h} \}_{h} \to E$ in $L^{1}$ and 
$\{ R_{h} \}_{h} \to 0^{+}$, 
we have to check that
$$\overline{F_0}(E) \leq \liminf_{n\to \infty} F_{R_{h}}(E_{h}).$$
We can suppose that
$$\liminf_{h \to \infty} F_{R_{h}}(E_{h})< +\infty$$
and we can
extract a subsequence $\{ E_{h_{k}} \}_{k}$ such that
$F_{R_{h_{k}} }(E_{h_{k}})$ is finite 
and
$$\liminf_{h \to \infty} F_{R_{h}}(E_{h}) = 
\lim_{k \to \infty} F_{R_{h_{k}}}(E_{h_{k}}).$$
Since $E_{h_{k}}$ belongs to ${\cal U}_{R_{h_{k}}}$, 
by Corollary 3.2 in \cite{bel-dm-pao}
it follows that
$$F_{R_{h_{k}}}(E_{h_{k}}) = \overline{F_0}(E_{h_{k}})$$
and then,
$$\overline{F_0}(E) \leq \lim_{k \to \infty} F_{R_{h_{k}}}
(E_{h_{k}});$$
hence, the liminf inequality is proved.

To obtain the limsup inequality, fixed $E\in {\cal M}$ 
and $\{ R_{h} \}_{h} \searrow 0$, we have 
to find a sequence $\{ E_{h} \}_{h} \to E$ in $L^{1}$ such that
$$\overline{F_0}(E) \geq  \limsup_{h \to \infty} F_{R_{h}}(E_{h}).$$
We can assume $\overline{F_0}(E)$ finite; then, 
there exists a sequence $\{ A_{k} \}_{k}$ such that
$A_{k}$ is in $C^{2}$, $\{ A_{k} \}_k$ converges to $E$ in $L^{1}$, 
as $k \to \infty$, and
$$\overline{F_0}(E) = \lim_{k \to \infty} F_0(A_{k}).$$
The smoothness of $A_{k}$ implies that there is $r_{k}>0$ such 
that $A_{k}$ belongs 
to the class ${\cal U}_{r_{k}}$.
Let us define by induction the following sequence of indices:
\begin{eqnarray*}
h_{1} & = & \min \{ h : R_{h} \leq r_{1}  \} \\
h_{k} & = & \min \{ h> h_{k-1} : R_{h} \leq r_{k} \}
\end{eqnarray*}
and the sets
$$
E_{h} := \begin{cases} 
A_{1} & \text{for $h< h_{1}$}, \\
A_{k} & \text{for $h_{k} \leq h < h_{k+1}$, $k\geq 1$}. 
\end{cases}$$
It is easy to verify that $\{ E_{h} \}_{h}$ is the required sequence.
\qed

\section{A variational problem in Image Segmentation} 

In this section we apply the results of the previous ones
to state an existence theorem for the {\sc Nitzberg} and {\sc Mumford}
problem in the class ${\cal U}_R$.
For every 
$k \in \enne$ and for every $E_{1},\ldots,E_{k}\in {\cal U}_{R}$ 
let us define the following functional:
\begin{multline}\label{bognone}
G_{k}(E_{1},\ldots,E_{k}) :=
 \alpha \int_{\Omega \setminus \cup_{i=1}^{k} E_{i}} 
|g - g_{\Omega 
\setminus \cup_{i=1}^{k} E_{i}}|^{2} dx + \\
+ \sum_{i=1}^{k}\left(
\alpha \int_{E'_{i} \cap \Omega} |g - 
g_{E'_{i} \cap \Omega}|^{2} dx 
  + \beta {\cal L}^{2}(E_{i}) + \gamma \int_{\partial E_{i}} \varphi(\kappa) 
d{\cal H}^{1}\right), 
\end{multline}
where $\alpha, \beta, \gamma$ are positive parameters,
$E'_i:=E_i\setminus \cup_{j=1}^{i-1}E_j$,  $g$ is a given function in 
$L^{2}(\Omega )$, $\varphi: \reali \to \reali$ is a given convex
function, and 
$\kappa$ denotes the curvature of 
$\partial E$.
If we take
$$\varphi(\kappa) = \begin{cases}
\nu + a \kappa^2 & \text{if $|\kappa |<\frac{b}{a}$}, \\
\nu + b |\kappa| & \text{if $|\kappa |\geq \frac{b}{a}$},
\end{cases}$$
we obtain exactly the original model proposed in \cite{nit-mum}.
 
\begin{theorem}\label{minimo1}
For every $R>0$ and for every  $k\in \enne$ the problem 
\begin{equation}\label{pb1}
\min \left\{ G_{k}(E_1, \ldots, E_k) : 
E_1, \ldots, E_k \in {\cal U}_{R} \right\}
\end{equation}
admits a solution.
\end{theorem}

{\sc Proof.} For the sake of simplicity, 
we perform the proof only for $k=1$; the general case 
follows by a similar argument, involving only some further 
difficulties of notation.

Let $\{ E_{m} \}_{m}$ be a minimizing sequence in ${\cal U}_{R}$ 
for the functional $G_{1}$. 
We can suppose that all non-empty connected components 
of each $E_{m}$ meet $\Omega$; indeed, if 
we call $\tilde{E}_m$ the union of the connected components of $E$ 
which intersect $\Omega$, we have
that $G_1 (\tilde{E}_m) \leq G_1 (E_m)$, 
and then, we can replace $E_m$ by $\tilde{E}_m$.
By Lemma \ref{diametri}
the sequence results equibounded.

Applying Theorems \ref{cpt} and \ref{lsc}
to the sequence $\{ E_{m}\}_m$, 
we obtain a subsequence $\{ E_{m_{h}} \}_{h}$ and a set 
$E \in {\cal U}_{R}$ such that
\begin{description}
\item[i)]  $E_{m_{h}} \to E$ in $L^{1}$ and a.e.;
\item[ii)]  $\displaystyle \int_{\partial E} \varphi(\kappa) 
d{\cal H}^{1} \leq \liminf_{h\to \infty} 
\int_{\partial E_{m_h}} \varphi(\kappa) 
d{\cal H}^{1}$.
\end{description}
We observe that
\begin{eqnarray*}
\lefteqn{\left| \int_{\Omega \cap E_{m_{h}}} 
|g - g_{\Omega \cap E_{m_{h}}}|^{2} dx 
- \int_{\Omega \cap E} |g - g_{\Omega \cap E}|^{2} dx
\right| \ \leq} \\
& \leq & \int_{\Omega} \chi_{\Omega \cap E_{m_{h}}} 
\left| |g - g_{\Omega \cap E_{m_{h}}}|^{2} 
- |g - g_{\Omega \cap E}|^{2} 
\right| dx 
+ \int_{\Omega} 
\left| \chi_{\Omega \cap E_{m_{h}}} -\chi_{\Omega \cap E} \right| 
|g - g_{\Omega \cap E}|^{2} dx;
\end{eqnarray*}
hence, applying the Dominated Convergence Theorem
to both addends, we can conclude that
$$\lim_{h\to \infty} \int_{\Omega \cap E_{m_{h}}} 
|g - g_{\Omega \cap E_{m_{h}}}|^{2} dx = 
\int_{\Omega \cap E} |g - g_{\Omega \cap E}|^{2} dx.$$
Analogously, 
$$\lim_{h \to \infty} \int_{\Omega \setminus E_{m_{k}}} 
|g - g_{\Omega \setminus E_{m_{h}}}|^{2} dx 
= \int_{\Omega \setminus E} |g - g_{\Omega \setminus E}|^{2} dx.$$
At this point it is clear that 
$$G_{1}(E) \leq \liminf_{h \to \infty} G_{1}(E_{m_{h}}),$$ 
and that $E$ minimizes the functional.
\qed

As explained in \cite{nit-mum-shi}, the integer $k$ is the number
of depth levels of the reconstructed image; denoting
by $(E_1, \ldots, E_k)$ the solution of (\ref{pb1}),
the set $E_i$ represents all the objects at the $i$-th level.
If $k$ is not a priori fixed, we can consider the variational problem
studied in the following theorem.

\begin{theorem}\label{minimo2}
For every $R>0$ the problem 
$$\min \left\{ G_{k}(E_1, \ldots, E_k) : 
E_1, \ldots, E_k \in {\cal U}_{R},\,  k\in \enne \right\}$$
admits a solution.
\end{theorem}

{\sc Proof.} Let $\{(E_{1}^{m},\ldots,E_{k_{m}}^{m})\}_m$ 
be a minimizing sequence.
Since for every $l\in\enne$,  
$j\in\{1,\ldots, l-1 \}$ and $A_{1},\ldots, A_{l-1} \in {\cal U}_{R}$
we have that
$$G_{l}(A_{1},\ldots, A_{j-1}, \emptyset, A_{j},\ldots, A_{l-1})=
G_{l-1}(A_{1},\ldots, A_{j-1}, A_{j},\ldots, A_{l-1}),$$
we can suppose that $E_{j}^{m}\neq\emptyset$ 
for every $m$ and for every $j\in\{1\ldots, k_m\}$,
and so,
$$G_{k_m}(E_{1}^{m},\ldots, E_{k_m}^{m})\geq \beta \sum_{j=1}^{k_m}
{\cal L}^{2}(E_{j}^{m})\geq \beta\, k_{m}{\cal L}^{2}(B(0,R));$$
therefore, the sequence $\{k_{m}\}_{m}$ must be bounded 
and so admits a constant subsequence: 
now we can conclude by applying Theorem \ref{minimo1}. \qed

If we are interested not only in detecting contours, but also
in cleaning and regularizing the image, we can
consider the following variational problem:
\begin{multline}\label{gradiente}
\min  \left\{ 
\alpha \int_{\Omega \setminus \cup_{i=1}^{k} \partial E'_{i}} |u-g|^{2}dx
+ \sum_{i=1}^{k}\left( \beta {\cal L}^{2}(E_i) + 
\gamma \int_{\partial E_{i}} \varphi(\kappa )\, d {\cal H}^{1}\right)  
+ \delta 
\int_{\Omega \setminus \cup_{i=1}^{k} \partial E'_{i}} 
|\nabla u|^{2} dx : \right. \\
\left.  \phantom{\sum_{0}^{0} \int}
u \in C^{1}(\Omega \setminus \cup_{i=1}^{k} \partial E'_{i}),\, 
E_{1},\dots, E_{k}
\in {\cal U}_{R} \right\}, 
\end{multline}
where $k$ is fixed in $\enne$, $\delta$ is a positive parameter 
and we use the same notation as before.

\begin{theorem}
Let $g$ be a function in $L^{\infty}(\Omega)$. 
Then, for every $R>0$ and for every
$k\in \enne$ the problem in (\ref{gradiente})
admits a solution.
\end{theorem}

{\sc Proof.} We first look for a solution $(u, E_{1},\dots, E_{k})$ where
$u \in W^{1,2}(\Omega \setminus \cup_{i=1}^{k} \partial E'_{i})$.

Let $\{(u_{h}, E_{1}^{h}, \dots, E_{k}^{h})\}_{h}$ be a minimizing sequence
for the functional in (\ref{gradiente}).
By a truncation argument we can suppose that $\| u_h \|_{\infty}
\leq \| g \|_{\infty}$ and, as in
the proof of Theorem \ref{minimo1}, we can assume
that
$\{ E_{i}^{h} \}_h$ is 
equibounded for every $i+1, \ldots, k$. By Theorem \ref{cpt}
there exist $E_{1}, \dots, E_{k}$ belonging to
${\cal U}_{R}$ such that, up to subsequences, 
$$ E_{i}^{h} \kconv E_{i}  \qquad \hbox{and}
\qquad E_{i}^{h}\to E_{i} \ \hbox{in} \ L^{1}(\Omega).$$
Arguing as in the proof of Theorem \ref{cpt}, one can easily check that
if $U$ is an open subset compactly contained in
$\Omega \setminus \cup_{i=1}^{k} \partial E'_{i}$, then
for $h$ large $U$ is compactly contained in
$\Omega \setminus \cup_{i=1}^{k} \partial (E_{i}^{h})'$.

Since $\{u_{h}\}_{h}$ is equibounded in $W^{1,2}(U)$, up to subsequences, 
there exists $u\in W^{1,2}(U)$ such that $u_{h} \rightharpoonup u$ in
w-$W^{1,2}(U)$.
By the weakly lower semicontinuity of the $L^{2}$-norm, by Theorems \ref{cpt}
and \ref{lsc},
we obtain
\begin{eqnarray*}
\lefteqn{\alpha \int_{U}|u-g|^{2}dx
+ \sum_{i=1}^{k}\left(
\beta {\cal L}^{2}(E_i)+
 \gamma
\int_{\partial E_{i}} \varphi(\kappa )\, d {\cal H}^{1}\right)
+ \delta \int_{U}|\nabla u|^{2}dx \ \leq}  \\
& \leq &  \liminf_{h \to \infty} \, 
 \alpha \int_{U}|u_{h}-g|^{2}dx
+\sum_{i=1}^{k}\left(
\beta {\cal L}^{2}(E_{i}^{h})+ \gamma \int_{\partial E_{i}^{h}} 
\varphi(\kappa )\, d {\cal H}^{1}\right)
+ \delta \int_{U} |\nabla u_{h}|^{2}dx \\
& \leq &  \liminf_{h \to \infty}  \,
\alpha \int_{\Omega \setminus \cup_{i=1}^{k} \partial (E_{i}^{h})'} 
|u_{h}-g|^{2}dx 
+ \sum_{i=1}^{k}
\left(
\beta {\cal L}^{2}(E_{i}^{h})+ 
\gamma \int_{\partial E_{i}^{h}} \varphi(\kappa )
\,d {\cal H}^{1}\right) \\
& & + \delta \int_{\Omega \setminus \cup_{i=1}^{k} \partial (E_{i}^{h})'}
|\nabla u_{h}|^{2}dx.
\end{eqnarray*}
Let us construct a sequence of open subsets compactly contained
in $\Omega \setminus \cup_{i=1}^{k} \partial E'_{i}$ 
and increasing to it; the previous argument combined with 
a diagonal procedure allows us to conclude that there exists $u\in W^{1,2}
(\Omega \setminus \cup_{i=1}^{k} \partial E'_{i})$ such that  
$(u, E_{1}, \dots,
E_{k})$ minimizes the functional.

Since $u-g \in L^{\infty}(\Omega \setminus 
\cup_{i=1}^{k} \partial E'_{i})$, the regularity theory
for elliptic equations ensures that
$u\in W^{2,p}_{loc}(\Omega \setminus 
\cup_{i=1}^{k} \partial E'_{i})$ for every $p< \infty$, hence
$u \in C^{1}(\Omega \setminus 
\cup_{i=1}^{k} \partial E'_{i})$.
\qed

Let us suppose now that $k$ is not a priori fixed: 
arguing as in Theorem \ref{minimo2}, we can prove
the following result.

\begin{theorem}
Let $g$ be a function in $L^{\infty}(\Omega)$. Then,
for every $R>0$ the problem
\begin{multline}
\min  \left\{ 
\alpha \int_{\Omega \setminus \cup_{i=1}^{k} \partial E'_{i}} |u-g|^{2}dx
+ \sum_{i=1}^{k}\left( \beta {\cal L}^{2}(E_i) + 
\gamma \int_{\partial E_{i}} \varphi(\kappa )\, d {\cal H}^{1}\right)  
+ \delta 
\int_{\Omega \setminus \cup_{i=1}^{k} \partial E'_{i}} 
|\nabla u|^{2} dx : \right. \\
\left. \phantom{\sum \int}
 u \in C^{1}(\Omega \setminus \cup_{i=1}^{k} \partial E'_{i}),\, 
E_{1},\dots, E_{k}
\in {\cal U}_{R}, \, k\in \enne \right\} 
\end{multline}
admits a solution.
\end{theorem}

We conclude this section by giving an example of 
non trivial (i.e. non empty) minimizer.

\begin{example}
Let us set $g:= \chi_{B(0,R)}$ and assume
$R>1$.
For a suitable choice of $\Omega$ and of the parameters 
$\alpha, \beta, \gamma$,
$B(0,R)$ minimizes the functional 
$$G_1 (E) = \alpha \int_{\Omega \setminus E} 
|g - g_{\Omega \setminus E}|^{2} dx
+ \alpha \int_{E \cap \Omega} |g - 
g_{E \cap \Omega}|^{2} \,dx 
  + \beta {\cal L}^{2}(E) + \gamma \int_{\partial E} (1 + |\kappa |^2)
\, d{\cal H}^{1}.$$
\end{example}

{\sc Proof.} It is known (see Theorem 5.7.3 in \cite{car}) that for 
every smooth closed curve $\gamma$, it results that
\begin{equation}\label{docarmo}
2\pi \leq \int_{\gamma} |\kappa| \, d{\cal H}^{1}.
\end{equation}
Holder inequality and (\ref{docarmo}) imply that
for every $E\in {\cal U}_{R}$, $E\neq \emptyset$, the following
inequality holds:
$$\int_{\partial E} |\kappa |^{2} d{\cal H}^{1} \geq 
\frac{4\pi^{2}}{{\cal H}^{1}(\partial E)},$$
so that 
$$G_{1}(E)\geq \beta \pi R^{2} + \gamma \left( {\cal H}^{1}(\partial E)
+ \frac{4\pi^{2}}{{\cal H}^{1}(\partial E)} \right).$$
Since ${\cal H}^{1}(\partial E) \geq 2\pi R$ and $R>1$,
$$G_{1}(E) \geq
\beta \pi R^{2} + \gamma \left( 2\pi R + \frac{2\pi}{ R} \right)
=  G_{1}(B(0,R)).$$
Finally,
$$G_{1}(\emptyset)=
\alpha \left( \pi R^{2} - \frac{\pi^{2}R^{4}}{{\cal L}^{2}(\Omega)}\right)
\geq \beta \pi R^{2} + \gamma \left( 2\pi R + \frac{4\pi^{2}}{2\pi R} \right)
$$
for a suitable choice of $\Omega$ and of the parameters.
\qed

\section*{Acknowledgements}
We wish to thank Professor Gianni Dal Maso for having suggested to us
the study of this problem and for fruitful
discussions.

\end{document}